\documentclass[MSNbibl,seceqn,nameyear,dvips]{arxbj}
\usepackage{mathrsfs}

\aid{0}
\volume{19}
\issue{4}
\pubyear{2013}
\firstpage{1243}
\lastpage{1267}
\doi{10.3150/12-BEJSP18} 

\makeatletter

\newcommand{\cH}{\mathcal H}
\newcommand{\cA}{\mathcal A}
\newcommand{\IT}{\mathbb T}

\newcommand{\cI}{\mathcal I}

\newcommand{\IP}{\mathbb P}

\newcommand{\IE}{\mathbb E}
\newcommand{\IR}{\mathbb R}

\newcommand{\cX}{\mathcal X}
\newcommand{\cY}{\mathcal Y}

\newcommand{\cS}{\mathcal S}

\newcommand{\IZ}{\mathbb Z}

\newtheorem{theorem}{Theorem}[section]
\newproclaim{definition}{Definition}[section]

\renewcommand{\phi}{\varphi}
\newcommand{\D}{\Delta}

\newcommand{\dtv}{d_{\mathrm{TV}}}

\newcommand{\bigo}{\mathrm{O}}
\newcommand{\lito}{\mathrm{o}}

\newcommand{\toinf}{\to\infty}
\newcommand{\tozero}{\to0}

\def\blambda{{\bolds\lambda}}

\newcommand{\Po}{\operatorname{Po}}
\newcommand{\CP}{\operatorname{CP}}
\newcommand{\Bi}{\operatorname{Bi}}
\newcommand{\N}{\mathrm{N}}

\newcommand{\I}{\mathrm{I}}

\newcommand{\Var}{\operatorname{Var}}
\newcommand{\law}{\mathscr{L}}

\def\bmid{\mid}

\makeatother

\begin{document}
\begin{frontmatter}

\title{Approximating dependent rare events}
\runtitle{Approximating dependent rare events}

\begin{aug}
\author[1]{\fnms{Louis H. Y.} \snm{Chen}\corref{}\thanksref{1}\ead[label=e1]{matchyl@nus.edu.sg}} \and
\author[2]{\fnms{Adrian} \snm{R\"ollin}\thanksref{2}\ead[label=e2]{adrian.roellin@nus.edu.sg}}
\runauthor{L.H.Y. Chen and A. R\"ollin} 
\address[1]{Department of Mathematics, National University of
Singapore, 10 Lower Kent Ridge Road, Singapore
119076.
\printead{e1}}
\address[2]{Department of Statistics and Applied Probability, National
University of Singapore, 6 Science Drive 2, Singapore
117546. \printead{e2}}
\end{aug}


%
\begin{abstract}
In this paper we give a historical account of the development of
Poisson approximation using Stein's method and present some of the main
results. We give two recent applications, one on maximal arithmetic
progressions and the other on bootstrap percolation. We also discuss
generalisations to compound Poisson approximation, Poisson process
approximation and multivariate Poisson approximation, and state a few
open problems.
\end{abstract}

%

\end{frontmatter}

\section{Introduction}\label{sec1}

The Poisson limit theorem as commonly found in textbooks of probability states
that the number of successes in $n$ independent trials converges in
distribution to a Poisson distribution with mean $\lambda> 0$ if the maximum
of the success probabilities tends to $0$ and their sum converges to
$\lambda$. The case where the trials have equal success probabilities was
implicitly proved by Abraham \citet{Moivre1712} in his solution to the problem
of finding the number of trials that gives an even chance of
getting $k$ successes. However, it was Sim\'eon-Denis \citet
{Poisson1837} who first gave an
explicit form of the Poisson distribution and proved the limit theorem
for independent trials with equal success probabilities, that is, for
the binomial
distribution. The Poisson distribution was not much used before Ladislaus
\citet{Bortkiewicz1898} expounded its mathematical properties and
statistical usefulness.

In his book \emph{Ars Conjectandi}, published posthumously in 1713, Jacob
Bernoulli (1654--1705) considered games of
chance and
urn models with two
possible outcomes and proved what is now known as the weak law of large
numbers. He stressed that the probability of winning a game or of
drawing a
ball of a particular color from an urn (with replacement) remains the same
when the game or the drawing of a ball is repeated. This has led to the use
of the term Bernoulli trials\vadjust{\eject} to represent independent trials with the same
probability of success. Representing success by $1$ and failure by $0$, a
random variable taking values $0$ and $1$ is called a Bernoulli random
variable. However, in this article, a set or a sequence of Bernoulli random
variables need not be independent nor take the value 1 with equal
probabilities.
Also, if the success probability of a Bernoulli random variable is
small, the
event corresponding to success is called \emph{rare}.

The Poisson limit theorem suggests that the distribution of a sum of
independent Bernoulli random variables with small success probabilities
can be
approximated by the Poisson distribution with the same mean if the success
probabilities are small and the number of random variables is large. A
measure of the accuracy of the approximation is the total variation distance.
For two distributions $P$ and $Q$ over $\IZ_+ = \{0,1,2,\ldots\}$, the total
variation distance between them is defined by
\[
\dtv(P,Q) = \sup_{A\subset\IZ_+}\bigl| P(A) - Q(A)\bigr|, 
\]
which is also equal to
\[
\frac12\sup_{| h| =1} \biggl|\int h\,dP - \int h\, dQ \biggr| =
\frac12\sum_{i\in\IZ} \bigl| P\{i\}- Q\{i\} \bigr|.
\]
For the binomial distribution $\Bi(n, p)$, \citet{Prohorov1953} proved that
\[
\dtv \bigl(\Bi(n,p), \Po(np) \bigr) \leq p \biggl[\frac{1}{\sqrt {2\pi e}} + \bigo
\biggl(1 \wedge\frac{1}{\sqrt{np}}+p \biggr) \biggr],
\]
where $\Po(np)$ denotes the Poisson distribution with mean $np$. Here,
following Barbour, Holst and Janson (\citeyear{Barbour1992}), the formulation corrects a minor
error in the original paper. This result is remarkable in that the
approximation is good so long as $p$ is
small, regardless of how large $np$ is.

The result of Prohorov was generalised by \citet{LeCam1960} to
sums of independent Bernoulli random variables $X_1, \ldots, X_n$ with success
probabilities
$p_1,\ldots, p_n$ that are not necessarily equal. Let $W = \sum X_i$ and
$\lambda= \sum p_i$. Using the method of convolution operators,
\citet{LeCam1960} obtained the error bounds
%
\begin{equation}
\label{1} \dtv \bigl(\law(W), \Po(\lambda) \bigr) \leq\sum
_{i=1}^n p_i^2,
\end{equation}
and
%
\begin{equation}
\label{2} \dtv \bigl(\law(W), \Po(\lambda) \bigr) \leq\frac{8}{\lambda
}\sum
_{i=1}^n p_i^2
\qquad\mbox{if }\max_{1\leq i\leq n} p_i \leq{
\frac{1}{4}}.
\end{equation}
In terms of order, the bound in (\ref{1}) is better than that
in (\ref{2}) if
$\lambda<1$ and vice versa if $\lambda\geq1$. Combining (\ref{1}) and
(\ref{2}),
one obtains
a bound of the order $(1\wedge\lambda^{-1})\sum p_i^2$, which is
small so long as $\max p_i$ is small, regardless of how large $\lambda
$ is.
This form of the error bound has become the characteristic of Poisson
approximation in subsequent developments of the subject.
\eject

In this article we will discuss the use of Stein's ideas in the Poisson
approximation to the distributions of sums of dependent Bernoulli random
variables, its historical development, applications, and some generalisations
and open problems. The article is not intended to be a survey paper but an
exposition with a focus on explaining Stein's ideas and presenting some
results and recent applications. The references are not exhaustive
but contain only those papers that are relevant to the objective of this
article.

This paper is organised as follows. Section \ref{sec2} is a brief introduction to
Stein's method. Section~\ref{sec3} gives a brief overview of two approaches to Poisson
approximation using Stein's method, and Sections \ref{sec4} and \ref{sec5} discuss the
developments of these two approaches. Section~\ref{sec6} is devoted to two recent
applications of Poisson approximation and Section \ref{sec7} discusses three
generalisations of Poisson approximation.

\section{Stein's method}\label{sec2}

In his seminal 1972 paper published in the Sixth Berkeley Symposium,
Charles Stein introduced a new method of normal approximation. The
method did
not involve Fourier analysis but hinged on the solution of a differential
equation. Although the method was developed for normal approximation, Stein's
ideas were very general and the method was modified by \citet{Chen1975} for
Poisson approximation.
Since then the method has been constantly developed and applied to many
approximations beyond normal and Poisson and in finite as well as infinite
dimensional spaces. It has been applied in many areas including computational
biology,
computer science, combinatorial probability, random matrices,
reliability and many more. The method, together with its applications,
continues to grow and remains a very active research area.
See, for example, \citet{Stein1986}, \citet{Arratia1990},
Barbour, Holst and Janson (\citeyear{Barbour1992}),
\citet{Diaconis2004}, Barbour and Chen (\citeyear{SteinsMethod2005a,SteinsMethod2005b}),
\citet{Chatterjee2005}, \citet{Chen2011}, \citet{Ross2011},
\citet{Shih2011}, \citet{Nourdin2012}.

In a nutshell, Stein's method can be described as follows. Let $W$
and $Z$ be
random elements taking values in a space $\cS$ and let $\cX$ and $\cY
$ be
some classes of real-valued functions defined on $\cS$. In approximating
the distribution $\law(W)$ of $W$ by the distribution $\law(Z)$
of $Z$, we
write $\IE h(W) - \IE h(Z) = \IE Lf_h(W)$ for a test function $h\in\cY$,
where $L$ is a linear operator (Stein operator) from $\cX$ into $\cY$ and $f_h
\in\cX$ a
solution of the equation
\[
Lf = h - \IE h(Z)\qquad\mbox{(Stein equation).}
\]
The error
$\IE{Lf_h(W)}$ can then be bounded by studying the solution $f_h$ and
exploiting the probabilistic properties of $W$. The operator $L$
characterises $\law(Z)$ in the sense that $\law(W) = \law(Z)$ if and
only if for a
sufficiently large class of functions $f$ we have
\[
\IE{Lf(W)}= 0 \qquad\mbox{(Stein identity).}
\]

In normal approximation, where $\law(Z)$ is the standard normal distribution,
the operator used by \citet{Stein1972} is given by $Lf(w) = f'(w) - wf(w)$
for $w\in\IR$, and in Poisson approximation, where $\law(Z)$ is the Poisson
distribution with mean $\lambda> 0$, the operator $L$ used by
\citet{Chen1975} is given by $Lf(w) = \lambda f(w + 1) - wf(w)$ for $w
\in
\IZ_+$. However the operator $L$ is not unique even for the
same approximating distribution but depends on the problem at hand. For
example, for normal approximation $L$ can also be taken to be the
generator of
the Ornstein--Uhlenbeck process, that is, $Lf(w) = f''(w) - wf'(w)$, and for
Poisson approximation, $L$ taken to be the generator of an immigration-death
process, that is, $Lf(w) = \lambda[f(w + 1) - f(w)] + w[f(w - 1) - f(w)]$.
This generator approach, which is due to \citet{Barbour1988}, allows
extensions to multivariate and process settings.
Indeed, for multivariate normal approximation, $Lf(w) = \D f(w) -
w\cdot
\nabla f(w)$, where $f$ is defined on the Euclidean space; see
\citet{Barbour1990} and \citet{Goetze1991}, and also \citet
{Reinert2009} and
\citet{Meckes2009}.

\section{Poisson approximation}\label{sec3}

In Poisson approximation, the main focus has been on bounding the total
variation distance between the distribution of a sum of dependent Bernoulli
random variables and the Poisson distribution with the same mean. One
of the
main objectives has been to obtain a bound which is the ``correct''
generalisation of the bound obtained by \citet{LeCam1960},
specifically, one
with the multiplicative factor $1\wedge\lambda^{-1}$.

Broadly speaking, there are two main approaches to Poisson
approximation using Stein's method, the local approach and the
size-bias coupling approach. The local approach was first studied by
\citet{Chen1975} and developed further by Arratia, Goldstein and
Gordon (\citeyear{Arratia1989,Arratia1990}), who presented Chen's
results in a form which is easy to use and applied them to a wide range
of problems including problems in extreme values, random graphs and
molecular biology. The size-bias coupling approach dates back to
\citet{Barbour1982} in his work on Poisson approximation for
random graphs. \citet{Barbour1992} presented a systematic
development of monotone couplings, and applied their results to random
graphs and many combinatorial problems. A~recent review of Poisson
approximation by \citet{Chatterjee2005} used Stein's method of
exchangeable pairs to study classical problems in combinatorial
probability. They also reviewed a size-bias coupling of Stein
[(\citeyear{Stein1986}), p. 93] for any set of dependent Bernoulli
random variables.

\section{The local approach}\label{sec4}

The operator $L$ given by $Lf(w) = f'(w) - wf(w)$ for $w \in\IR$,
which was
used by \citet{Stein1972} for normal approximation, is constructed by showing
that $\IE\{f'(Z) - Zf(Z)\} = 0$ for all bounded absolutely continuous
functions $f$ if $Z \sim\N(0,1)$.
This identity is proved by integration by parts. As a discrete counterpart,
the operator $L$ given by $Lf(w) = \lambda f(w + 1) - wf(w)$ for $w \in
\IZ_+$, which was used by \citet{Chen1975} for Poisson approximation, is
constructed by showing
that $\IE\{\lambda f(Z) - Zf(Z)\} = 0$ for all bounded real-valued
functions $f$ if $Z\sim\Po(\lambda)$, using summation by parts.

Using the Stein equation
%
\begin{equation}
\label{3} \lambda f(w + 1) - wf(w) = h(w) - \IE h(Z),
\end{equation}
where $|h| = 1$ and $Z$ has the Poisson distribution with mean $\lambda
> 0$,
\citet{Chen1975}
developed Stein's method for Poisson approximation for sums of
$\phi$-mixing sequences of Bernoulli random variables $X_1,\ldots,X_n$ with
success probabilities $p_1,\ldots,p_n$. When specialised to independent
Bernoulli random variables, his results yield
\[
\dtv \bigl(\law(W),\Po(\lambda) \bigr) \leq3 \biggl(1\wedge\frac1{\sqrt{
\lambda}} \biggr) \sum_{i=1}^n
p_i^2
\]
and
\[
\dtv \bigl(\law(W),\Po(\lambda) \bigr) \leq\frac{5}{\lambda}\sum
_{i=1}^n p_i^2,
\]
where $W=\sum X_i$. These results improve slightly those of
\citet{LeCam1960}.

Chen's proofs depend crucially on the bounds he obtained on the
solution of
(\ref{3}) and its smoothness. These bounds were improved by
\citet{Barbour1983}, who proved that for $h = I_A$, $A \subset
\IZ_+$,
%
\begin{equation}
\label{4}
\Vert f_h\Vert _\infty\leq1\wedge\frac{1.4}{\sqrt{\lambda}}
\end{equation}
and
%
\begin{equation}
\label{5} \Vert\D f_h\Vert _\infty \leq\frac{1 - e^{-\lambda}}{\lambda}
\leq1\wedge\frac{1}{\lambda},
\end{equation}
where \mbox{$\Vert\cdot\Vert _\infty$} denotes the supremum norm and $\D
f(w) = f(w + 1)
- f(w)$.

It is perhaps instructive to see how easily Le Cam's results, with
substantially
smaller constants, can be proved by Stein's method using (\ref{5}).

Let $W$ be the sum of independent Bernoulli random variables $X_1,
\ldots, X_n$ with success probabilities $p_1, \ldots, p_n$, and let
$W^{(i)} = W - X_i$ for $i = 1, \ldots, n$. For any bounded
real-valued function~$f$,
\begin{eqnarray*}
\IE\bigl\{\lambda f(W + 1) - Wf(W)\bigr\} & = & \sum_{i=1}^n
\IE \bigl\{p_i f(W+1) - X_if(W) \bigr\}
\\
& = &\sum_{i=1}^n p_i \IE
\bigl\{ f(W+1) - f \bigl(W^{(i)}+1 \bigr) \bigr\}
\\
& = &\sum_{i=1}^n p_i \IE
\bigl\{ X_i \Delta f \bigl(W^{(i)}+1 \bigr) \bigr\}
\\
& = &\sum_{i=1}^n p_i^2
\IE\D f \bigl(W^{(i)}+1 \bigr).
\end{eqnarray*}
By choosing $f = f_h$, a bounded solution of (\ref{3}), where $h =
I_A$ and $A
\subset
\IZ_+$, we
obtain
%
\begin{eqnarray}
\label{6} \dtv \bigl(\law(W),\Po(\lambda) \bigr) & = & \sup_{A\subset\IZ_+}
\bigl|\IP[W \in A] - \IP[Z \in A]\bigr|
\nonumber\\[-8pt]\\[-8pt]
& \leq &\Vert\D f_h\Vert _\infty\sum
_{i=1}^n p_i^2 \leq
\biggl(1\wedge\frac1\lambda \biggr) \sum_{i=1}^n
p_i^2.
\nonumber
\end{eqnarray}
We wish to remark that the solution $f_h$ is unique except at $w = 0$,
but the
value of $f_h$ at $w = 0$ is never used in the calculation. So it has been
conveniently set to be 0.

The above proof of (\ref{6}) is given in \citet{Barbour1984}, who
also proved that
\[
\dtv \bigl(\law(W),\Po(\lambda) \bigr) \geq\frac{1}{32} \biggl(1\wedge
\frac1\lambda \biggr) \sum_{i=1}^n
p_i^2.
\]
This shows that $(1\wedge\lambda^{-1}) \sum p_i^2$ is of the best
possible order for the Poisson approximation. Indeed, it has been
proved by
\citet{Deheuvels1986}, using a semigroup approach, and also by \citet{Chen1992}
and \citet{Barbour1995}, using Stein's method, that $\dtv (\law(W),
\Po(\lambda) )$ is asymptotic to $\sum p_i^2$
(respectively $(2\pi e)^{-1/2}
\lambda^{-1}\sum p_i^2$) as $\max p_i \tozero$ and $\lambda\tozero$
(respectively
$\lambda\toinf$).

We end this section by stating a theorem of
Arratia, Goldstein and
Gordon [(\citeyear{Arratia1989,Arratia1990}), Theorem 1], which was proved using
(\ref{4}) and
(\ref{5}).

%
%

\begin{theorem}\label{thm1} Let $\{X_\alpha\dvtx  \alpha\in J\}$ be Bernoulli
random variables with success probabilities $p_\alpha$, $\alpha\in
J$. Let
$W=\sum_{\alpha\in J} X_\alpha$ and
$\lambda= \IE W = \sum_{\alpha\in J} p_\alpha$. Then, for any
collection of
sets
$B_\alpha\subset J$, $\alpha\in J$,
\[
\dtv \bigl(\law(W),\Po(\lambda) \bigr) \leq \biggl(1\wedge\frac{1}{\lambda}
\biggr) (b_1+b_2 ) + \biggl(1\wedge\frac{1.4}{\sqrt{\lambda}}
\biggr) b_3
\]
and
\[
\bigl|\IP[W=0] - e^{-\lambda} \bigr| \leq \biggl(1\wedge \frac{1}{\lambda}
\biggr) (b_1+b_2+b_3),
\]
where
\begin{eqnarray*}
b_1 &=& \sum_{\alpha\in J}\sum
_{\beta\in B_\alpha} p_\alpha p_\beta, \qquad
b_2 = \sum_{\alpha\in J}\sum
_{\beta\in B_\alpha\setminus\{\alpha
\}} \IE(X_\alpha X_\beta),
\\
b_3 &=& \sum_{\alpha\in J} \bigl|
\IE(X_\alpha| X_\beta, \beta \notin B_\alpha)-
p_\alpha \bigr|.
\end{eqnarray*}
\end{theorem}

If for each $\alpha\in J$, $X_\alpha$ is independent of $\{X_\beta\dvtx
\beta\notin B_\alpha\}$, then $b_3 = 0$,
and we call $\{X_\alpha\dvtx  \alpha\in J\}$ locally dependent with
dependence neighbourhoods $\{B_\alpha\dvtx \alpha\in J \}$. An
$m$-dependent sequence of random variables, which is a special case of a
$\phi$-mixing sequence, is locally
dependent.

The wide applicability of Theorem \ref{thm1} is illustrated through many
examples in
Arratia, Goldstein and Gordon
(\citeyear{Arratia1989,Arratia1990}). Many problems
to which Theorem \ref{thm1} is applied are concerned with locally dependent
random variables.

\section{The size-bias coupling approach}\label{sec5}

In his monograph, Stein [(\citeyear{Stein1986}), pp. 89--93] considered the following
general problem of Poisson approximation. Let $X_1, \ldots, X_n$ be
dependent Bernoulli random variables with success probabilities $p_i =
\IP[X_i = 1]$ for $i = 1,\ldots,n$. Let $W = \sum X_i$ and let
$\lambda= \IE W$ with $\lambda> 0$. Assume $I$ to be uniformly
distributed over $\{1, \ldots,n\}$ and independent of $X_1, \ldots, X_n$.
Then for any bounded real-valued function $f$ defined on
$\{0,1,\ldots,n\}$,
%
\begin{equation}
\label{7} \IE\bigl\{Wf(W)\bigr\} = \lambda \IE \bigl(f(W) \bmid X_I
= 1 \bigr).
\end{equation}
If $W^*$ and $W$ are defined on the same probability space such that the
distribution of $W^*$ equals the conditional distribution of $W$
given $X_I =
1$, then (\ref{7}) becomes
\[
\IE\bigl\{Wf(W)\bigr\} = \lambda \IE f\bigl(W^*\bigr),
\]
from which one obtains
%
\begin{equation}
\label{8} \dtv \bigl(\law(W), \Po(\lambda) \bigr) \leq\bigl(1 - e^{-\lambda}
\bigr) \IE\bigl|W + 1 - W^*\bigr|.
\end{equation}
From (\ref{8}), one can see that if the distribution of $W +1$ is close
to that
of $W^*$, then the distribution of W is approximately Poisson with mean
$\lambda$, and (\ref{8}) gives an upper bound on the total variation distance.

This approach to Poisson approximation was reviewed in
Chatterjee, Diaconis and
  Meckes (\citeyear{Chatterjee2005}),
who also applied (\ref{8}) to a variety of problems, such as
the matching problem, the coupon-collector's problem and the birthday problem.

In their monograph, \citet{Barbour1992} studied Poisson
approximation for Bernoulli random variables satisfying monotone coupling
assumptions. We state their main theorem in this context as follows.

\begin{theorem}\label{thm2} Let $\{X_\alpha\dvtx  \alpha\in J\}$ be Bernoulli
random variables with success probabilities $p_\alpha$, $\alpha\in J$.
Suppose for each $\alpha\in J$, there exists $\{Y_{\beta,\alpha}\dvtx
\beta\in
J\}$ defined on the same probability space as $\{X_\alpha\dvtx  \alpha\in
J\}$ such that
\[
\law \bigl( \{Y_{\beta,\alpha}\dvtx  \beta\in J \} \bigr) = \law \bigl(
\{X_\alpha\dvtx  \alpha\in J \bmid X_\alpha= 1 \} \bigr).
\]
Let $W = \sum X_\alpha$, $\lambda= \IE W = \sum p_\alpha$, and
$Z \sim\Po(\lambda)$.
\begin{enumerate}
\item If
%
\begin{equation}
\label{9} Y_{\beta,\alpha} \leq X_\beta \qquad\mbox{for all $\beta\in
J$ (negatively related),}
\end{equation}
then
%
\begin{equation}
\label{10} \dtv \bigl(\law(W),\Po(\lambda) \bigr) \leq (1\wedge\lambda ) \biggl(
1 -\frac{\Var(W)}{\lambda} \biggr).
\end{equation}
\item If
%
\begin{equation}
\label{11} Y_{\beta,\alpha} \geq X_\beta \qquad\mbox{for all $\beta\in
J$ (positively related),}
\end{equation}
then
%
\begin{equation}
\label{12} \dtv \bigl(\law(W),\Po(\lambda) \bigr) \leq (1 \wedge \lambda )
\biggl( \frac{\Var(W)}{\lambda} -1 + \frac{2}{\lambda}\sum
_{\alpha\in J} p_\alpha^2 \biggr).
\end{equation}
\end{enumerate}
\end{theorem}

From (\ref{10}) and (\ref{12}), one can see that $\law(W)$ is approximately
$\Po(\lambda)$ if $\Var(W)/\lambda$ is close to~$1$.

The proof of Theorem \ref{thm2} is pretty similar to that for (\ref
{8}). Let
$V_\alpha= \sum_{\beta\neq\alpha} Y_{\beta,\alpha}$
and $W^{(\alpha)} = W -
X_\alpha$ for $\alpha\in J$. Then for any bounded real-valued function
$f$ defined on $\{0, 1, \ldots,| J| \}$,
%
\begin{eqnarray}
\label{13} \IE\bigl\{ Wf(W) \bigr\} & = & \sum_{\alpha\in J}
p_\alpha \IE \bigl(f\bigl(W^{(\alpha)} + 1\bigr) \bmid
X_\alpha= 1 \bigr)
\nonumber\\[-8pt]\\[-8pt]
& = &\sum_{\alpha\in J} p_\alpha \IE
f(V_\alpha+ 1) = \lambda \IE f(V_I + 1),
\nonumber
\end{eqnarray}
where $I$ is independent of all the $X_\alpha$ and $V_\alpha$, and
 $\IP[I =
\alpha] = p_\alpha/\lambda$, $\alpha\in J$.

Using the monotone properties (\ref{9}) and (\ref{11}), one gets
\[
\dtv \bigl(\law(W),\Po(\lambda) \bigr) \leq \bigl(1 - e^{-\lambda
} \bigr)
\IE\bigl[(W + 1) - (V_I + 1)\bigr]
\]
for the negatively related case, and
\[
\dtv \bigl(\law(W),\Po(\lambda) \bigr) \leq \bigl(1 - e^{-\lambda
} \bigr)
\bigl(\IE X_I + \IE\bigl[(V_I + 1) -
\bigl(W^{(I)} + 1\bigr)\bigr] \bigr)
\]
for the positively related case. Straightforward calculations then yield
(\ref{10}) and (\ref{12}).

\citet{Barbour1992} also established conditions for existence of
monotone couplings and applied Theorem \ref{thm2} to large number of problems
in random permutations, random graphs, occupancy and urn models,
spacings, and
exceedances and extremes.

The coupling approach of \citet{Stein1986} and of \citet{Barbour1992}
can actually be formulated under the general framework of size-bias coupling.
Here is the definition of size-biased distribution; see \citet{Goldstein1996}.

\begin{definition}\label{def1} Let $W$ be a non-negative random
variable with
mean $\lambda> 0$. We say that $W^s$ has the \emph{$W$-size biased
distribution} if
\[
\IE\bigl\{Wf(W)\bigr\} = \lambda \IE f\bigl(W^s\bigr)
\]
for all real-valued functions $f$ such that the expectations exist.
\end{definition}
\vspace*{-12pt}
\eject

If $W$ is a non-negative integer-valued random variable, then $\IP[W^s
= k] =
k\IP[W = k]/\lambda$ for $k \geq1$.
The following theorem follows immediately.

\begin{theorem}\label{thm3} Let $W$ be a non-negative integer valued random
variable with $\IE W = \lambda> 0$. Assume that $W^s$ and $W$ are
defined on
the same probability space, that is, assume that there is size-bias coupling.
Then we have
\[
\dtv \bigl(\law(W),\Po(\lambda) \bigr) \leq \bigl(1 - e^{-\lambda
} \bigr)\IE
\bigl|W + 1 - W^s\bigr|.
\]
\end{theorem}

Note that in the case where $W$ is a sum of Bernoulli random variables, $W^s$
can be taken to be $W^*$ in (\ref{8}) or $V_I+1$ in (\ref{13}).
Furthermore, it is clear from Theorem \ref{thm3} that the Poisson
distribution is the only distribution such that its size-biased distribution
is the original distribution shifted by one.

We conclude by saying that a large portion of the literature on the coupling
approach to Poisson approximation falls under the general framework of
size-bias coupling. Indeed, (\ref{13}) provides a general way for
constructing size-bias coupling for sums of Bernoulli random
variables. Couplings involving the size-biased distribution, however, have
found applications beyond Poisson approximation; see for example
\citet{Pekoz2011} and \citet{Pekoz2013}.

%
%

\section{Applications}\label{sec6}

A remarkable feature of Theorem \ref{thm1} and Theorem \ref{thm2} is
that the
error bounds depend only on the first two moments of the random
variables. It
also happens that many interesting scientific problems can be
formulated as
occurrences of dependent rare events. For example, one is often
interested in
the
maximum of a set of random variables $\xi_1,\ldots,\xi_n$. For a threshold
$t$, define $X_i = I[\xi_i>t]$ for $i=1,\ldots,n$, and let $W=\sum X_i$. Then
%
\begin{equation}
\label{14} P[\max\xi_i\leq t ] = \IP[W=0].
\end{equation}
Often $t$ is large, so that $\{ \xi_1>t\},\ldots,\{\xi_n>t\}$ are
rare events.
If the $X_1,\ldots,X_n$ satisfy the conditions of Theorem \ref{thm1} or
Theorem \ref{thm2} and the error bound is small, then
\[
P[\max\xi_i \leq t] \approx e^{-\lambda_t} \qquad\mbox{where }
\lambda_t = \sum_{i=1}^n P[
\xi_i>t].
\]
Since the appearance of Theorems \ref{thm1} and \ref{thm2}, Poisson
approximation has been applied to a large number of problems in many
different fields, which include computational biology, random graphs and
large-scale networks, computer science, statistical physics, epidemiology,
reliability theory, game theory, and financial mathematics. In computational
biology, Poisson approximation is typically used to
calculate $p$-values in
sequence comparison, while in random graphs, it is used to count the
copies of a small graph in a large graph. Here is a sample of
publications on
problems in different fields, in which Poisson approximation is applied:
\citet{Dembo1994},
\citet{Neuhauser1994},
\citet{Waterman1994},
\citet{Waterman1995},
\citet{Embrechts1997},
\citet{Karlin2000},
\citet{Barbour2001b},
\citet{Lange2002},
Lippert, Huang and
  Waterman (\citeyear{Lippert2002}),
\citet{Grimmett2003},
\citet{Franceschetti2006},
\citet{Hart2008},
\citet{Draief2010},
Falk, H{\"u}sler and Reiss (\citeyear{Falk2011}).

In what follows, we will present two recent applications of Poisson
approximation, one by \citet{Benjamini2007} on maximal arithmetic
progressions, and the other by \citet{Bollobas2012a} on bootstrap percolation.

\subsection{Maximal arithmetic progressions}\label{sec6.1}

The occurrences of arithmetic progressions in subsets of the set of positive
integers are of interest in number theory. \citet{Tao2007} gave a historic
account of the topic, in particular, Szemer\'edi's theorem, which
states that
any ``dense'' subset of positive integers must contain arbitrarily long
arithmetic progressions.

Benjamini, Yadin and Zeitouni (\citeyear{Benjamini2007,Benjamini2012}) analyse the
following probabilistic variant of arithmetic progressions. Let
$\xi_1,\ldots,\xi_n$ be i.i.d. Bernoulli random variables with success
probability $0<p<1$. We say
that there is an arithmetic progression of length at least $t$,
starting at
$a+s$ with a common difference $s$, if $\xi_a = 0$ and $\xi_{a+s}=
\xi_{a+2s}
= \cdots= \xi_{a+ts} = 1$ as long as $a+ts \leq n$.
Let $U_n$ denote the length of the maximal arithmetic progression among
$\xi_1,\ldots,\xi_n$. We have the following result.

\begin{theorem}[{[\citet{Benjamini2007}]}]\label{thm5} Let $x\in\IR$ be
fixed and
let $0\leq\delta_n<1$ be such that
 $x - \frac{2\log n}{\log p} + \frac{\log\log n}{\log p} - \delta
_n$ is
integer valued. Then
%
\begin{equation}
\label{15} \IP \biggl[U_n + \frac{2\log n}{\log p} - \frac{\log\log n}{\log
p} <
x - \delta_n \biggr] \sim\exp \biggl(\frac{(1-p)p^{x-\delta_n}\log p
}{4} \biggr)
\end{equation}
as $n\toinf$.
\end{theorem}

Note that the distribution of $U_n$ is of Gumbel-type. However, the
rounding effect of $\delta_n$ does not vanish, since
$U_n$ is integer-valued and, as one can show, the variance of $U_n$ is of
order $1$. Therefore, limiting distributions only exist along subsequences
$n_1, n_2, \ldots$
for which $\lim_{m\toinf}\delta_{n_m}$ exists, in which case the limiting
distribution is a discretised Gumbel distribution.

\begin{pf*}{Idea of proof}
Denote by $\cI_{n,t}$ the set of pairs $(a,s)$ of
positive integers that satisfy $a+ts\leq n$, and for each such
$(a,s)\in\cI_{n,t}$ define
\[
X_{a,s} = \I[\xi_a=0, \xi_{a+s} =
\xi_{a+2s} = \cdots= \xi _{a+ts} = 1].
\]
%
Let $W_{n,t} = \sum_{(a,s)} X_{a,s}$, were the sum ranges over all pairs
$(a,s)\in\cI_{n,t}$. Then $W_{n,t}$ counts the arithmetic
progressions of length at least t in $\{1, 2, \ldots, n\}$. Following
(\ref{14}),
we have
\[
\IP[U_n < t] = \IP[W_{n,t} = 0].
\]
We claim that
%
\begin{equation}
\label{16} \IP[W_{n,t} = 0] \approx e^{-\lambda_{n,t}},
\end{equation}
where
\[
\lambda_{n,t} = |\cI_{n,t}| qp^{t}
\]
with $q=1-p$. It is not difficult to see that
\[
|\cI_{n,t}| = \sum_{s=1}^{{\lfloor\frac{n-1}{t}\rfloor}}
(n-ts) \sim \frac{n^2}{2t}
\]
if $n,t\toinf$ as long as $t = \lito(n)$. We let
%
\begin{equation}
\label{17} t = x -\frac{2\log n}{\log p} + \frac{\log\log n}{\log p} - \delta_n,
\end{equation}
which is integer-valued by definition of $\delta_n$. Since with this
choice of
$t$ we have
\[
\lambda_{n,t} \sim\frac{-qp^{x-\delta_n}\log p }{4}
\]
as $n\toinf$, we have established (\ref{15}).

It remains to justify (\ref{16}) for $t$ defined as in (\ref{17}),
which we will
accomplish via Theorem \ref{thm1}. To this end,
let $A_{a,s} = \{ a,a+s,\ldots,a+ts \}$ for each $(a,s)\in\cI_{n,t}$.
Note that
$X_{a,s}$ and $X_{a',s'}$ are independent whenever the sets $A_{a,s}$ and
$A_{a',s'}$ are disjoint. Denote by $D_{a,s}(k)$ be number of pairs
$(a',s')\in\cI_{n,t}$ with $s'\neq s$, such that $| A_{a,s}\cap
A_{a',s'}|
= k$. From Benjamini, Yadin and Zeitouni
[(\citeyear{Benjamini2007}), Proposition 4] we have the estimate
\[
D_{a,s}(k) \leq %
\cases{ (t+1)^2 n, & if $k=1$,
\cr
(t+1)^2 t^2, & if $2\leq k\leq t/2+1$,
\cr
0, &
if $k>t/2+1$.} %
\]
We can now apply Theorem \ref{thm1}. Let
$N_{a,s}\subset\cI_{n,t}$ be the set of pairs $(a',s')$ such
that $A_{a,s}\cap
A_{a',s'}$ is non-empty. It is clear then that $b_3 = 0$. Now,
\begin{eqnarray*}
b_1 \leq\sum_{a,s} \Biggl(1+\sum
_{k=1}^{t} D_{a,s}(k) \Biggr)
p^{2(t+1)} = \bigo \biggl(\frac{n^2}{t}\cdot \bigl(1+t^2n
+t^5 \bigr)\cdot p^{2t} \biggr) = \lito(1),
\end{eqnarray*}
where we used that $p^{2t}=\bigo (\log(n) / n^4 )$. Since
$\IE(X_{a,s}X_{a',s'})= p^{2(t+1)-k}$ if $| A_{a,s}\cap
A_{a',s'}| =
k$, we also have
\[
b_2 \leq\sum_{a,s}\sum
_{k=1}^{t} D_{a,s}(k) p^{2(t+1)-k} =
\bigo \biggl(\frac{n^2}{t}\cdot \bigl(t^2n+t^4
p^{-t/2} \bigr)\cdot p^{2t} \biggr) = \lito(1).
\]
Hence, by Theorem \ref{thm1}, $|\IP[W_{n,t}=0]-e^{-\lambda
_{n,t}}| \tozero$
as $n\toinf$, justifying (\ref{16}).
\end{pf*}

We refer to Benjamini, Yadin and Zeitouni (\citeyear{Benjamini2007,Benjamini2012}) and \citet{Zhao2012}
for further
details and refinements.

\subsection{Bootstrap percolation}\label{sec6.2}

Consider the $d$-dimensional torus lattice $\IT_{n}^d = \IZ^d /
n\IZ^d$, along with the canonical $\ell_1$ distance, that is, the
smallest number of edges connecting to points. Two sites are connected
if their $\ell_1$ distance is $1$. \citet{Bollobas2012a}
considered \emph{$d$-neighbour bootstrap percolation} on $\IT_n^d$,
a~special type of a \emph{cellular automaton}. A vertex can be either
\emph{infected} or \emph{uninfected}. At each time step, an uninfected
vertex becomes infected if $d$ or more of its neighbours are infected
(at each time step, this rule is applied \emph{simultaneously} for all
vertices). Once a vertex is infected, it stays infected.

The rules of cellular automata are usually deterministic, and in the model
considered by \citet{Bollobas2012a}, randomness is added only at the
beginning: at time $0$, each vertex is
infected with probability $p$ and remains uninfected with probability
$q=1-p$, independently of all other vertices. With $A_t\subset
\IT_n^d$
denoting the set of all infected sites at time $t$,
we shall be interested in the first time
\[
T_n = \inf\bigl\{t\geq0: A_t = \IT_n^d
\bigr\}
\]
when all the sites are infected.

The following result says that, if we let $p$ converge to $1$ at
the right speed as $n\toinf$, $T_n$ is essentially concentrated on one
or two points.
In order to formulate the result, we define the combinatorial quantity
\[
m_t = \sum_{r=0}^t\sum
_{j=0}^r \pmatrix{d\cr j}.
\]
%
\begin{theorem}[{[Bollob{\'a}s \textit{et~al.} (\citeyear{Bollobas2012a}), Theorem 3]}] Fix a positive integer $t$.
If, for some function
$\omega(n)\toinf$,
%
\begin{equation}
\label{18} \biggl(\frac{\omega(n)}{n^{d}} \biggr)^{\frac{1}{m_{t-1}}} \leq q_n
\leq \biggl(\frac{1}{\omega(n)n^{d}} \biggr)^{\frac{1}{m_t}}
\end{equation}
then $\IP _{p_n}[T_n=t] \to1$. If instead, for some slowly varying function
$\omega(n)$,
%
\begin{equation}
\label{19} \biggl(\frac{1}{\omega(n) n^{d}} \biggr)^{\frac{1}{m_{t}}} \leq q_n
\leq \biggl(\frac{\omega(n)}{n^{d}} \biggr)^{\frac{1}{m_t}}
\end{equation}
then $\IP _{p_n}[T_n\in\{t,t+1\}]\to1$.
\end{theorem}

\begin{pf*}{Idea of proof}
As in the previous application, we reformulate the
problem as an extremal problem. Although we assume
that $t$ is fixed throughout, the arguments can
be extended to $t=\lito(\log n / \log\log n)$.
For each
$i\in\IT_n^d$, let
$Y_i$ be the time when vertex $i$ becomes infected, that is
\[
Y_i = \inf\{t\geq0: i\in A_t\}.
\]
Now clearly $T_n = \max_{i\in\IT_n^d} Y_i$. For each $i\in\IT
_n^d$, let
$X_{t,i} = \I[Y_i>t]$ be the indicator that vertex $i$ is uninfected
at time
$t$. Note that, although $Y_i$ and $Y_j$ are not independent for
any $i$ and
$j$, the indicators $X_{t,i}$ and $X_{t,j}$ are independent whenever the
$\ell_1$-distance between $i$ and $j$ is larger than $2t+1$, since infections
can only propagate an $\ell_1$-distance $1$ per time step.
With $W_{n,t}= \sum_{i\in\IT_n^d} X_{t,i}$, we have
%
\begin{equation}
\IP _p[T_n\leq t] = \IP _p[W_{n,t}=0].
\end{equation}
We claim that
\begin{equation}
\label{20}
\IP _p[W_{n,t}=0] \approx e^{-\lambda_{n,t}},
\end{equation}
where
\[
\lambda_{n,t} = \sum_{i\in\IT_n^d}\IP
_p[\mbox{$i$ is uninfected at time $t$}] = n^d
\rho_{n,t}(p)
\]
with $\rho_{n,t}(p) = \IP _p[\mbox{0 is uninfected at time $t$}]$.

Bollob{\'a}s \textit{et~al.} [(\citeyear{Bollobas2012a}),
Theorem 17] gave the following results about the behaviour of $\rho_{n,t}(p)$. If there
exists $C = C(t,d) >0$ such that
%
\begin{equation}
\label{21} q_n^{m_t} \leq\frac{C}{n^d},
\end{equation}
for all $n$, then
%
\begin{equation}
\label{22} \rho_{n,t}(p_n) \sim d^3
2^{d-1} q_n^{m_t}
\end{equation}
as $n\toinf$. Hence, if
%
\begin{equation}
\label{23}
q_n^{m_t} \leq\frac{1}{n^d \omega(n)}
\end{equation}
for some function $\omega(n)\toinf$, we have that $\lambda
_{n,t}(p_n) = n^d
\rho_{n,t}(p_n)\tozero$, so that, under (\ref{23}),
%
\begin{equation}
\label{24} \IP _{p_n}[T_n\leq t] \to1.
\end{equation}
If, in contrast, we have
%
\begin{equation}
\label{25}
q_n^{m_t} \geq\frac{\omega(n)}{n^d}
\end{equation}
for some function $\omega(n)\toinf$, we can argue as follows. A simple
coupling argument yields that
the system is monotone, that is, if $\tilde{p}\leq p$, we have $\rho
_{n,t}(\tilde{p})
\geq\rho_{n,t}(p)$, and hence $\lambda_{n,t}(\tilde{p}) \geq
\lambda_{n,t}(p)$.
Since by (\ref{22}) we have
\[
\lambda_{n,t}(p_n)\sim C d^3
2^{d-1}
\]
for arbitrarily large $C$, we must have
\[
\lambda_{n,t} (p_n )\toinf
\]
if (\ref{25}) is true, thus yielding
%
\begin{equation}
\label{26} \IP _{p_n}[T_n\leq t] \to0.
\end{equation}
Since the first inequality in (\ref{18}) is just (\ref{25}) with $t$
replaced by $t-1$, we have from (\ref{26}) that $\IP _{p_n}[T_n\leq
t-1] \to
0$. On the other hand, the second inequality of (\ref{18}) is just
(\ref{23}), hence
(\ref{24}) implies $\IP _{p_n}[T_n\leq t]\to1$. Thus, $\IP
_{p_n}[T_n=t] \to
1$. The proof of the second statement is analogous by observing that
(\ref{19})
implies
\[
\biggl(\frac{\tilde{\omega}(n)}{n^{d}} \biggr)^{\frac{1}{m_{t-1}}} \leq q_n \leq
\biggl(\frac{1}{\tilde{\omega}(n)n^{d}} \biggr)^{\frac{1}{m_{t+1}}},
\]
where, with $\alpha= m_{t-1}/m_t < 1$,
\[
\tilde{\omega}(n) = \frac{n^{d(1-\alpha)}}{\omega(n)^\alpha} \toinf.
\]

It remains to justify (\ref{20}). Again, by monotonicity it is enough to
consider (\ref{21}), since $\tilde{p} \leq p$ implies $\IP _{\tilde
{p}}[T_n\leq
t]\leq
\IP _{ p}[T_n\leq t]$. Let
\[
\tilde{\rho}_{n,t}(p) = \max_{j: d(0,j)\leq2t } \IP
_p[\mbox{$0$ and $j$ are uninfected at time $t$}].
\]
Bollob{\'a}s \textit{et al.} [(\citeyear{Bollobas2012a}), Lemma 19]
showed that, if (\ref{21}) holds, then
\[
\tilde{\rho}_{n,t}(p_n) = \lito \bigl(\rho_{n,t}(p_n)
\bigr).
\]
Let now $N_i = \{ j\in\IT_n^d:  d(i,j)\leq2t\}$. It is clear that
$X_{t,i}$ is independent of $(X_{t,j})_{j\notin N_i}$, hence $b_3 =
0$. With
the crude bound $| N_i| \leq t^d$, we have
\[
\frac{b_1}{\lambda_{n,t}(p_n)} \leq\frac{n^d t^d
\rho_{n,t}(p_n)^2}{\lambda_{n,t}(p_n)} = t^d
\rho_{n,t}(p_n) = \lito(1)
\]
and
\[
\frac{b_2}{\lambda_{n,t}(p_n)} \leq\frac{n^d t^d
\tilde{\rho}_{n,t}(p_n)}{\lambda_{n,t}(p_n)} = t^d\frac{\tilde{\rho}_{n,t}(p_n)}{\rho_{n,t}(p_n)} =
\lito(1),
\]
justifying (\ref{20}).
\end{pf*}

\section{Generalisations and open problems}\label{sec7}

In this section we will discuss three generalisations of Poisson approximation
and touch briefly on two other generalisations, the three generalisations
being compound Poisson approximation, Poisson process approximation and
multivariate Poisson approximation. Compound Poisson distributions on the
real line, the distributions of Poisson point processes, and multivariate
Poisson distributions are all compound Poisson distributions if viewed
in an
appropriate way, but the three approximations have been studied separately
because of the different contexts in which they arise and the different
problems to which they are applied.

\subsection{Compound Poisson approximation}\label{sec7.1}

In many probability models (see \citet{Aldous1989}), events occur in
clumps at
widely
scattered localities or at long irregular intervals in time. In
such situations, the Poisson approximation for the number of events occurring
either fails or performs poorly. If the number of clumps is approximately
Poisson, the clumps are roughly independent and their sizes close to
identically distributed, then the number of events occurring can be
approximated by a compound Poisson distribution. A typical example of events
occurring in clumps is earthquakes exceeding certain magnitude. Often
such an
earthquake is followed by a quick succession of several earthquakes before
normalcy is resumed.

We illustrate further the notion of clumps by presenting the example of the
longest head run discussed in Arratia, Goldstein and Gordon
(\citeyear{Arratia1989,Arratia1990}).
Note that this example is a special case of the maximal arithmetic
progressions in Section \ref{sec6.1}. Suppose a coin is tossed repeatedly where the
probability of falling heads is $p$ ($0 < p < 1$). Let $R_n$ be the
length of
the longest run of heads starting from within the first $n$ tosses.
What is
the asymptotic distribution of $R_n$ as $n\toinf$?

Let $Z_1, Z_2, \ldots$ be independent Bernoulli random variables with success
probability $p$ ($0 < p < 1$), where $\{Z_i = 1\}$ represents the event that
the coin falls heads at the $i$th toss. Let $J = \{1, 2, \ldots, n\}$
and let
$t
\geq1$. Define $Y_i = Z_{i}Z_{i+1}\cdots Z_{i+t-1}$ for $i = 1, 2,
\ldots, n$,
and define
\[
X_i = %
\cases{ Y_1, & if $i = 1$,
\cr
(1 -
Z_{i-1})Y_i, & if $2 \leq i \leq n$.} %
\]
Let $W = \sum X_i$ and let $\lambda= \IE W$. Then $\{R_n < t\} = \{W =
0\}$.

Define $B_i = \{ j \in J:  |i - j| \leq t\}$, $i = 1, 2, \ldots, n$.
Then $\{X_i: i \in J\}$ is locally dependent with dependence neighbourhoods
$\{B_i: i \in J\}$. Applying Theorem \ref{thm1}, we obtain $b_3 = b_2
= 0$,
and $b_1 < \lambda^2(2t+1)/n + \lambda p^t$.

Hence
%
\begin{equation}
\label{27} \bigl|\IP[R_n < t] - e^{- \lambda} \bigr| \leq \biggl(1
\wedge\frac
{1}{\lambda} \biggr) \bigl(\lambda^2(2t+1)/n + \lambda
p^t \bigr).
\end{equation}
Requiring that $\lambda$ remains bounded away from $0$ and from
$\infty$ and that the error bound tends to $0$ as $n\toinf$ leads to
the following conclusion: for a fixed integer $c$, $\IP[R_n -
{\lfloor\log _{1/p}(n(1-p))\rfloor} < c] \to\exp\{- p^{c-r}\}$ along a
subsequence of $n$ if and only if $\log_{1/p}(n(1-p)) -
{\lfloor\log_{1/p}(n(1-p))\rfloor} \to r \in[0,1]$ along the same
subsequence.\vspace*{1pt}

Now let $V = \sum Y_i$ and let $\mu= \IE V$. Then we also have $\{
R_n <
t\} =
\{V = 0\}$. The difference between the $X_i$ and the $Y_i$ is that while $X_i$
indicates a run of at least $t$ heads starting from the $i$th toss
preceded by
a
tail, $Y_i$ indicates a run of at least $t$ heads starting from
the $i$th toss
regardless of what precedes it. For a run of more than $t$ heads starting
from
the $i$th toss, say, $Z_{i-1} = 0$, $Z_{i} = \cdots= Z_{i+m-1}= 1$,
$Z_{i+m}=0$, where $m > t$, $X_i = 1$, $X_{i+1} = \cdots= X_{i+m-t} = 0$,
whereas $Y_i = Y_{i+1} = \cdots= Y_ {i+m-t } = 1$. Thus while $W$
counts the
clumps, which consist of runs of at least $t$ heads, $V$ counts the
clumps and
their sizes. The way the $X_i$ are defined so that $W$ counts only the clumps
is
called
\emph{declumping}.

If we apply Theorem \ref{thm1} to $V$, we will obtain a bound on $|\IP
[R_n <
t] - e^{-
\mu}|$. Since $\{Y_i: i \in J\}$ is locally dependent, $b_3 = 0$. But
$b_2$ does
not tend to $0$ if we require $\mu$ to be bounded away from $0$ and from
$\infty$.
Thus Poisson approximation fails. However, \citet{Arratia1990} showed
that the
distribution of $V$ is approximately compound Poisson
through an extension of Poisson approximation to Poisson process
approximation.

We pause for a moment to remark that there are two equivalent representations
of the compound Poisson distribution on $\IZ_+ = \{0, 1, 2, \ldots\}
$. Let
$\xi_1, \xi_2, \ldots$ be i.i.d. positive integer-valued random variables
with $\IP[\xi_1 = k] = \gamma_k$ for $k = 1, 2, \ldots\,$, and let $N$
be a Poisson
random variable with mean $\nu> 0$, independent of the $\xi_i$. The
distribution of $\xi_1 + \xi_2 + \cdots+ \xi_N$, which is compound Poisson,
is the same as that of $\sum iZ_i$, where the $Z_i$ are independent Poisson
random variables with means $\nu\gamma_i$ respectively. Let $\gamma$
be the
common distribution of the $\xi_i$. Then $\gamma= \sum\gamma
_i\delta_i$, where $\delta_i$ is the Dirac measure at $i$. We denote
this compound Poisson
distribution by $\CP(\nu\gamma) = \CP (\sum\nu\gamma
_i\delta_i )$ and call $\nu\gamma$ the generating measure.

\citet{Arratia1990} showed that by representing $\{Y_i: i \in J\}$ as a
Bernoulli process indexed by $J \times\{1, 2, \ldots\}$ where $J$ denotes
the location of clumps and $\{1, 2, \ldots\}$ the clump sizes, $\{Y_i:
i \in
J\}$ can be approximated in total variation by a Poisson process, which is
a collection of independent Poisson random variables indexed by $J
\times\{1,
2, \ldots\}$. By taking an appropriate projection and using the above
alternative representation of the compound Poisson distribution,
\citet{Arratia1990} obtained a bound on the total variation distance between
the distribution of V and a compound Poisson distribution. This in turn
provides an error bound for $|\IP[R_n < t] - e^{- \nu}|$, where $\nu
$ is the
mean of the Poisson number of terms in the compound Poisson
distribution and
is less than $\mu= \IE V$. This error bound is of the same order as
that in
(\ref{27}), but without the factor $1 \wedge\lambda^{-1}$. However,
it leads
to the same asymptotic distribution for $R_n$ as $n \toinf$
because $\lambda$
is bounded away from $0$ and from $\infty$.

The factor $ 1\wedge\lambda^{-1}$ is lost because Poisson process
approximation for the Bernoulli process representing $\{Y_i: i \in J\}$
requires too much information extraneous to the compound Poisson approximation
for $V$. A direct approach using Stein's method, which partially
recovers the
factor $1\wedge\lambda^{-1}$, was developed by \citet{Barbour1992c}. Let
$\lambda_i \geq0$, $i = 1, 2, \ldots$ such that $\sum\lambda_i <
\infty$. \citet{Barbour1992c} used the Stein equation
%
\begin{equation}
\label{28} \sum i \lambda_i f(w + i) - w f(w) = \I(w \in A)
- \IP[Z \in A]\qquad \mbox{for $w \in\IZ_+$,}
\end{equation}
where $A$ is a subset of $\IZ_+$, $Z_i$, $i \geq1$, are independent
$\Po(\lambda_i)$, and $Z
= \sum iZ_i$.

By solving (\ref{28}) analytically as well as writing $f(w) = g(w) - g(w
-1)$ and
using the generator approach to solve (\ref{28}), they obtained the
following bounds on the solution~$f_A$. For $A \subset\IZ_+ = \{0, 1,
2, \ldots\}$,
%
\begin{equation}
\label{29} \Vert f_A\Vert _{\infty}\leq \bigl(1\wedge
\lambda_1^{-1} \bigr)e^{\nu}, \qquad \Vert\D
f_A\Vert _{\infty}\leq \bigl(1\wedge\lambda_1^{-1}
\bigr)e^{\nu},
\end{equation}
where $\nu= \sum\lambda_i$, and if $i\lambda_i \downarrow0$, then
\[
\Vert f_A\Vert _{\infty}\leq %
\cases{ 1,& if $
\lambda_1 - 2\lambda_2 \leq1$,
\cr
\displaystyle
\frac{2}{(\lambda_1 - 2\lambda_2)^{1/2}} - \frac
{1}{\lambda_1 - 2\lambda_2},& if $\lambda_1 - 2
\lambda_2 > 1$,} %
\]
and
\[
\Vert\D f_A\Vert _{\infty} \leq1\wedge \biggl(
\frac{1}{4(\lambda-
2\lambda_2)^2}+ \frac{ \log^+ (2(\lambda_1 - 2\lambda
_2) )}{\lambda_1 - 2\lambda_2} \biggr).
\]

As in the case of Poisson approximation, the solution $f_A$ is unique
except at $w = 0$. Since its value at $w = 0$ is never used in the
calculation, it has been conveniently set to be $0$. Using the bounds
on $\|\D f_A\|_{\infty}$, \citet{Barbour1992c}
proved the
following theorem for locally dependent Bernoulli random variables.

\begin{theorem}\label{thm6}
Suppose $\{X_\alpha: \alpha\in J\}$ are locally dependent
Bernoulli random variables with success probabilities $p_\alpha$ and
dependence neighbourhoods $B_\alpha\subset C_\alpha$, $\alpha\in
J$, such
that for each $\alpha\in J$, $X_\alpha$ is independent of $\{
X_\beta:
\beta\in B_\alpha^c\}$ and $\{X_\beta: \beta\in B_\alpha\}$ is
independent of $\{X_\beta: \beta\in C_\alpha^c\}$. Let $ W = \sum X_\alpha$ and let $Y_\alpha= \sum_{\beta\in B_\alpha} X_\beta$.
Define $\lambda_i = i^{-1}\sum\IE X_\alpha\I[Y_\alpha= i]$ for $i
= 1, 2,
\ldots\,$, let $\nu= \sum\IE X_\alpha Y_\alpha^{-1} = \sum\lambda
_i$, and let $\gamma= \sum(\lambda_i/\nu)\delta_i$.
\begin{enumerate}
\item We have
\[
\dtv \bigl(\law(W), \CP(\nu\gamma) \bigr) \leq \bigl(1 \wedge
\lambda_1^{-1} \bigr) e^ \nu\sum
_{\alpha\in J}\sum_{\beta
\in
C_\alpha} p_\alpha
p_\beta.
\]
\item If $i\lambda_i \downarrow0$ as $i \to\infty$, then we have
\begin{eqnarray*}
&&\dtv \bigl(\law(W), \CP(\nu\gamma) \bigr)
\\
&&\quad\leq2 \biggl[1\wedge \biggl(\frac{1}{4(\lambda_1 -
2\lambda_2)^2} + \frac{\log^+ (2(\lambda_1 - 2\lambda_2) )}{\lambda_1 -
2\lambda_2}
\biggr) \biggr]\sum_{\alpha\in J}\sum
_{\beta\in C_\alpha} p_\alpha p_\beta.
\end{eqnarray*}
\end{enumerate}
\end{theorem}

If $\lambda_i = 0$ for $i \geq3$ and $\lambda_1 < 2\lambda_2$, it can
be shown that both $\Vert f_A\Vert _{\infty}$ and $\Vert\D f_A\Vert
_{\infty}$ grow exponentially fast with $\nu$ (see \citet
{Barbour1998a,Barbour1999a}). This shows that the bounds in (\ref{29})
cannot be much improved. To circumvent this difficulty Barbour and Utev
(\citeyear{Barbour1998a,Barbour1999a}) considered bounds on
\begin{eqnarray*}
H_0^a(\nu\gamma) &:=& \sup_{A \subset\IZ_+}
\sup_{w > a}\bigl|f_A(w)\bigr|,
\\
H_1^a(\nu\gamma) &:=& \sup_{A \subset\IZ_+}
\sup_{w > a}\bigl|f_A(w+1)- f_A(w)\bigr|.
\end{eqnarray*}
Assuming that the generating function of $\gamma= \sum(\lambda_i/\nu
)\delta_i$ has a radius of convergence $R > 1$ and assuming some other
conditions,\vadjust{\eject} \citet{Barbour1999a} proved that there exist
constants $C_0$, $C_1$ and $C_2$ depending on $\gamma$ such that for
any $a > C_2\nu m_1 + 1$, where $m_1$ is the mean of $\gamma$,
\[
H_0^a(\nu\gamma) \leq C_0
\nu^{-1/2},
\]
and
\[
H_1^a(\nu\gamma) \leq C_1
\nu^{-1}.
\]

The expressions for $C_0, C_1$ and $C_2$ are complicated but explicit.
Sufficient conditions can be found under which these constants are
uniformly bounded. Using the bound on $H_1^a(\nu\gamma)$, \citet
{Barbour2000c} proved the following theorem.

\begin{theorem}
For $n \geq1$, let $\lambda_{in}> 0$ for $i = 1, 2, \ldots\,$.
Let $W_n$, $n = 1, 2, \ldots\,$, be a sequence of non-negative,
integer-valued random variables such that for each $n \geq1$ and each
bounded $f:\IZ_+ \to\IR$,
\[
\biggl|\IE \biggl(\sum_{i\geq1}i\lambda_{in}f(W_n
+ i) - W_nf(W_n) \biggr) \biggr| \leq\Vert\D f \Vert
_{\infty
}\epsilon_n.
\]
Let $\nu_n = \sum_{i \geq1}\lambda_{in} < \infty$ and $\gamma
_{in} = \lambda_{in}/\nu_n$. Assume that
\begin{eqnarray*}
&&\mbox{\textup{\hphantom{ii}(i)}}\quad \lim_{n \to\infty}\gamma_{in} =
\gamma_i\mbox{ for each }i \geq1,\qquad
\hspace*{29.7pt}\mbox{\textup{(ii)}}\quad \inf_{n\geq1}\gamma_{1n}>0,
\\
&&\mbox{\textup{(iii)}}\quad \sup_{n \geq1}\sum
_{i \geq1} \gamma _{in}R^i < \infty\mbox{ for
some }R > 1,
\qquad
\mbox{\textup{(iv)}}\quad \inf_{n \geq1} \nu_n > 2.
\end{eqnarray*}
Then there exist positive constants $K < \infty$ and $c < 1$ such that
for any $x$ satisfying $c < x < 1$ and any $n$ for which\/ $\IE W_n
\geq(x - c)^{-1}$,
\[
\dtv \bigl(\law(W_n), \CP(\Gamma_n) \bigr) \leq K(1 -
x)^{-1} \bigl(\nu_n^{-1}\epsilon_n + \IP
\bigl(W_n \leq(1+x)\IE W_n/2 \bigr)\bigr),
\]
where the generating measure $\Gamma_n = \sum\lambda_{in}\delta_i =
\nu_n \sum\gamma_{in}\delta_i$.
\end{theorem}
In their efforts to obtain bounds on the solution of the Stein
equation (\ref{28}) so that the bounds resemble or ``correctly''
generalise those in the Poisson approximation, Barbour and Xia
[(\citeyear{Barbour1999}), Theorem
2.5] obtained the following theorem by treating compound
Poisson approximation as a perturbation of Poisson approximation.

\begin{theorem}\label{thm7} Let $\lambda_i \geq0$, $i \geq1$, satisfy
\[
\theta:= \frac{1}\lambda\sum i(i-1) \lambda_i < {
\frac{1}{2}} \qquad \mbox{where }\lambda= \sum i
\lambda_i < \infty.
\]
Then for any subset $A\subset\IZ_+$, the Stein equation (\ref{28}) has
a bounded
solution $f = f_A$ satisfying
\[
\Vert f_A\Vert _\infty\leq\frac{1}{(1 - 2\theta)
\lambda^{1/2}}, \qquad \Vert
\Delta f_A\Vert _\infty\leq\frac{1}{(1 - 2\theta)
\lambda}.
\]
\end{theorem}

Using the bound on $\Vert\Delta f_A\Vert _\infty$ for the locally dependent
Bernoulli random variables defined in Theorem \ref{thm6}, we obtain the
following theorem.

\begin{theorem}\label{thm8} Let $\{X_\alpha: \alpha\in J\}$ be locally
dependent Bernoulli random variables as defined in Theorem \ref{thm6}. Let $W
= \sum X_\alpha$ and let $Y_\alpha= \sum_{\beta\in B_\alpha}
X_\beta$.
Define $\lambda_i = i^{-1}\sum\IE X_\alpha\I[Y_\alpha= i]$ for $i
= 1, 2, \ldots\,$. If $\theta: = \lambda^{-1}\sum i(i-1) \lambda
_i < {\frac{1}{2}}$, where $\lambda= \sum i\lambda_i < \infty$,
and $\gamma=\sum (\lambda_i/\nu)\delta_i$, then
\[
\dtv \bigl(\law(W), \CP(\nu\gamma) \bigr) \leq 2 \biggl(1\wedge
\frac{1}{(1 - 2\theta) \lambda} \biggr)\sum_{\alpha\in
J} \sum
_{\beta\in
C_\alpha} p_\alpha p_\beta.
\]
\end{theorem}

Much progress has been made on bounding the solution of the Stein
equation (\ref{28}) in compound Poisson approximation. The results
presented in this section are quite satisfactory although many
conditions on the $\lambda_i$ or the generating measure are required.
It still remains a tantalising question as to what general results one
could obtain by using a different Stein equation or by using a
non-uniform bound on its solution, and to what extent one could do away
with those conditions on the $\lambda_i$. \citet{Roos2003c} used the generating
function approach of Kerstan to study compound Poisson approximation
for sums
of independent random variables without imposing any condition on
the $\lambda
_i$, but the method works only under the condition of independence.
Even for sums of independent random variables it is unclear if the
results of \citet{Roos2003c} can be proved using Stein's method. For
further reading on compound Poisson approximation, see \citet
{Barbour2001c}, and \citet{Erhardsson2005}.

\subsection{Poisson process approximation}\label{sec7.2}

In Poisson process approximation, both the number of rare events that occur
and the respective locations at which they occur are approximated by a Poisson
point process on a metric space. In the longest head run example
discussed in
\citet{Arratia1990}, the information on the locations
where the events occur is used in the calculation of the compound Poisson
approximation. In \citet{Leung2005}, a Poisson process
approximation for palindromes in a DNA is used to provide a
mathematical basis
for modelling the palindromes as i.i.d. uniform random variables on an
interval. The total variation distance is used for the Poisson process
approximation in the longest head run example, but in general such a distance
is not appropriate. For example, the total variation distance between a
Bernoulli process indexed by $\{i/n: i = 1, 2, \ldots, n\}$ with success
probability $\lambda/n$ and a Poisson process on $[0, 1]$ with
rate $\lambda
$ is always $1$, although the Bernoulli process converges weakly to the
Poisson process as $n \to\infty$. A distance which is commonly used in
process approximations is the Wasserstein distance.

By writing $f(w) = g(w) - g(w - 1)$, \citet{Barbour1988} converted the
Stein equation (\ref{3}) to a second order difference equation and
introduced the
generator approach to extend Poisson approximation to higher dimensions
and to
Poisson process approximation. Following the generator approach,
\citet{Barbour1992b} established a general framework for Poisson
process approximation. In this framework, a compact metric
space $\Gamma$
endowed with a metric $d_0 \leq1$ is the carrier space, $\Xi$ is a point
process on $\Gamma$ with finite intensity measure $\blambda$ of total mass
$\lambda$, where $\blambda(A) = \IE\Xi(A)$ for every Borel set
in $\Gamma$,
and $Z$ is a Poisson point process on $\Gamma$ with the same intensity measure
$\blambda$. Let $\cX$ be the configuration space $\{\sum_{1 \leq i
\leq k}
\delta_{\alpha_i}: \alpha_i \in\Gamma, k \geq0\}$. Define a metric $d_1
\leq1$ on $\cX$ by
\[
d_1 \biggl( \sum_{1\leq i\leq m}
\delta_{x_i}, \sum_{1\leq i\leq n}
\delta_{y_i} \biggr) = %
\cases{ 1, & if $m \neq n$,
\cr
\displaystyle n^{-1}\min\sum d_0 (x_i,
y_{\pi(i)} ), & if $m = n$,} %
\]
where the minimum is taken over all permutations $\pi$ of $\{1, 2,
\ldots
n\}$.

Define $\cH= \{h:\cX\to\IR: |h(\xi_1) - h(\xi_2)| \leq d_1(\xi_1,
\xi_2)\}$. The Wasserstein distance with respect to $d_1$ between
the distributions of two point processes $\Xi_1$ and $\Xi_2$
on $\Gamma$
with finite intensity measures is defined by
\[
d_2 \bigl(\law(\Xi_1), \law(\Xi_2) \bigr) =
\sup_{h\in\cH} \bigl|\IE h(\Xi _1) - \IE h(
\Xi_2)\bigr|.
\]
Note that $d_2$ is a metric bounded by $1$.
The Stein equation for approximating the distribution of the point process
$\Xi$ by that of the Poisson point process $Z$ is
%
\begin{eqnarray}
\label{30} \cA g(\xi) &:=& \int_\Gamma\bigl[g(\xi+
\delta_x) - g(\xi)\bigr]\blambda (dx) + \int\bigl[g(\xi-
\delta_x) - g(\xi)\bigr]\xi(dx)
\nonumber\\[-8pt]\\[-8pt]
& = & h(\xi) - \IE h(Z),
\nonumber
\end{eqnarray}
where $h \in\cH$ and $\cA$ is the generator of a measure-valued
immigration-death process $Y_\xi(t)$ with immigration
intensity $\blambda$,
per capita unit death rate, $Y_\xi(0) = \xi$, and stationary distribution
$\law(Z)$.

The Stein equation (\ref{30}) has a solution $g = g_h$ given by
\[
g_h(\xi) = - \int_{0}^\infty\bigl[Eh
\bigl(Y_\xi(t)\bigr) - Eh(Z)\bigr]\,dt.
\]
Using coupling, \citet{Barbour1992b} obtained the following bounds
on $g_h$:
%
\begin{eqnarray}
\label{31} \bigl|\Delta_\alpha g_h(\xi)\bigr| &:=& \bigl|g_h(
\xi+ \delta_\alpha) - g_h(\xi)\bigr| \leq1\wedge1.65
\lambda^{-1/2},
\\
\bigl|\Delta^2_{\alpha\beta} g_h(\xi)\bigr| &:=&
\bigl|g_h(\xi+ \delta_\alpha +\delta _\beta) -
g_h(\xi+ \delta_\alpha) - g_h(\xi+
\delta_\beta) + g_h(\xi)\bigr|
\nonumber\\[-8pt]\\[-8pt]
&\leq& 1\wedge\frac{5(1 +
2\log^+(2\lambda/5)}{2\lambda},\nonumber
\end{eqnarray}
where $\lambda$ is the total mass of $\blambda$.

In applications, the logarithmic term in (\ref{31}) carries over to the
error bounds in the approximation. Attempts were made to remove the
logarithmic terms. Xia (\citeyear{Xia1997,Xia2000}) succeeded in some
special cases. A general result in the form of a non-uniform bound on
$|\Delta_{\alpha\beta} ^2g_h(\xi )|$ was obtained by
\citet{Brown2000} and later improved by \citet{Xia2005},
which is given as
%
\begin{equation}
\label{32} \bigl|\Delta_{\alpha\beta} ^2g_h(\xi)\bigr| \leq1
\wedge \biggl(\frac{3.5}{\lambda} + \frac{2.5}{|\xi| + 1} \biggr),
\end{equation}
where $|\xi|$ is the number of points in $\xi$, that is, the total measure
of $\xi$.

Using (\ref{32}), the error bound on the Wasserstein distance for Poisson
process approximation for Bernoulli processes has the same factor as
that on
the total variation distance for the Poisson approximation for sums of
independent Bernoulli random variables, namely, $1\wedge\lambda^{-1}$.

\citet{Chen2004} studied Stein's method for Poisson process
approximation from the point of view of Palm theory. For a point process
$\Xi$ on $\Gamma$ with finite intensity measure, the Palm
process $\Xi
_\alpha$ associated with $\Xi$ at $\alpha\in\Gamma$ has the same
distribution as the conditional distribution of $\Xi$ given that a
point has
occurred at $\alpha$. A point process $\Xi$ on $\Gamma$ with finite
intensity measure $\blambda$ is locally dependent with neighbourhoods
$\{A_\alpha: \alpha\in\Gamma\}$ if $\law (\Xi_\alpha
^{(\alpha)} )
= \law (\Xi^{(\alpha)} )$ $\blambda$-a.s., where $\Xi
_\alpha
^{(\alpha)} $ and $\Xi^{(\alpha)}$ are respectively the restrictions
of $\Xi
_\alpha$ and $\Xi$ to $A_\alpha^c$ for each $\alpha\in\Gamma$.

The following theorem, which uses (\ref{32}), is Corollary 3.6 in
\citet{Chen2004}.

\begin{theorem}
Let $\Xi$ be a locally dependent point process on the compact
metric space $\Gamma$ with finite intensity measure $\blambda$ and with
neighbourhoods $ \{A_\alpha: \alpha\in\Gamma\}$, and let $Z$ be a
Poisson point process on $\Gamma$ with the same intensity
measure $\blambda$.
Let $\lambda$ be the total measure of $\blambda$. Then
\begin{eqnarray*}
d_2 \bigl(\law(\Xi), \law(Z) \bigr) & \leq & \IE\int
_{\alpha\in\Gamma} \biggl(\frac{5}{\lambda} + \frac{3}{|\Xi^{(\alpha)}| + 1} \biggr)
\bigl(\Xi (A_\alpha)-1\bigr)\Xi(d\alpha)
\\
&&{} + \int_{\alpha\in\Gamma} \int_{\beta\in A_\alpha} \biggl(
\frac{5}{\lambda} + \IE\frac{3}{|\Xi^{(\alpha\beta)}| +
1} \biggr)\blambda(d\alpha)\blambda(d
\beta),
\end{eqnarray*}
where $|\xi|$ is the total measure of $\xi$ and $\xi^{(\alpha\beta
)}$ is
the restriction of $\xi$ to $A_\alpha^c \cap A_\beta^c$.
\end{theorem}

This theorem gives the factor $ 1\wedge\lambda^{-1}$ in the Wasserstein
distance error bound for the Poisson approximation for Bernoulli
Processes. It
has also been applied to Poisson process approximation for palindromes
in a
DNA in \citet{Leung2005}, and to Poisson point process
approximation for the Mat\'ern hard-core process in \citet{Chen2004}.

For further reading on Poisson process approximation, see \citet{Xia2005}.

\subsection{Multivariate Poisson approximation}\label{sec7.3}

For the multivariate analogue of Poisson approximation, we consider
independent Bernoulli random d-vectors, $X_1, \ldots, X_n$ with
\[
\IP \bigl[X_j = e^{(i)} \bigr] = p_{j,i}, \qquad
\IP[X_j = 0] = 1 - p_j, \qquad 1 \leq i \leq d, 1 \leq j
\leq n,
\]
where $e^{(i)}$ denotes the $i$th coordinate vector in $\IR^d$
and $p_j =
\sum_{1 \leq i \leq d} p_{j,i}$.
\eject

Let $W = \sum X_j$, $\lambda= \sum p_j$, $\mu_i = \lambda
^{-1}\sum_{1
\leq j \leq n} p_{j,i}$, and let $Z = (Z_1, \ldots, Z_d)$,
where $Z_1, \ldots, Z_d$ are independent Poisson random variables with means
$\lambda\mu_1, \ldots, \lambda\mu_d$. Using the Stein equation
%
\begin{eqnarray}
\label{33} \cA g(j) &=& \sum\lambda\mu_i \bigl\{g \bigl(j +
e^{(i)} \bigr) - g(j) \bigr\} + \sum j^{(i)} \bigl\{g
\bigl(j - e^{(i)} \bigr) - g(j) \bigr\}
\nonumber\\[-8pt]\\[-8pt]
& = & \I[j \in A] - \IP[Z \in A],
\nonumber
\end{eqnarray}
where $A$ is a subset of $\IZ_+^d$ and $\cA$ the generator of a multivariate
immigration-death process whose stationary distribution is
$\law(Z)$, \citet{Barbour1988} proved that
%
\begin{equation}
\label{34} \dtv \bigl(\law(W), \law(Z) \bigr) \leq\sum
_{1 \leq j \leq n} p_j^2\wedge \biggl(
\frac{c_\lambda}{\lambda} \sum_{1 \leq i \leq d} \frac{ p_{j,i}^2}{\mu_i}
\biggr),
\end{equation}
where $c_\lambda= {\frac{1}{2}}+ \log^+(2\lambda)$.

The error bound in (\ref{34}) looks like the ``correct''
generalisation of
$(1\wedge\lambda^{-1})\sum_{1 \leq j \leq n} p_j^2$ in the univariate
case except for the factor $c_\lambda$, which grows logarithmically with
$\lambda$.

Using the multivariate adaption of Kerstan's generating function method,
\citet{Roos1999} proved that
%
\begin{equation}
\label{35} \dtv \bigl(\law(W), \law(Z) \bigr) \leq8.8 \sum
_{1 \leq j \leq n} p_j^2\wedge \biggl(
\frac{1}{\lambda} \sum_{1 \leq i \leq d} \frac{ p_{j,i}^2}{\mu_i}
\biggr),
\end{equation}
which improves over (\ref{34}) in removing $c_\lambda$ from the error bound
although the absolute constant is increased to 8.8.

The error bound in (\ref{34}) was obtained by bounding $\Delta_{ik}g_A$
in the
error term in the approximation where $g_A$ is the solution of the Stein
equation (\ref{33}), $\Delta_ig(k) = g (k + e^{(i)} ) -
g(k)$ and
$\Delta_{ik} =
\Delta_i(\Delta_k)$. By studying the behaviour of $\Delta_{ik}g_A$,
\citet{Barbour2005} showed that the order of the bound in (\ref{34})
is best
possible for $d \geq2$ if it is proved by bounding $\Delta_{ik} g_A$.
By an indirect approach to bounding the error term \citet{Barbour2005} obtained
two error
bounds, one of which comes very close to (\ref{35}) and the other
better than
an earlier bound of \citet{Roos1998}.

There does not seem to be much progress on multivariate Poisson approximation
using Stein's method since 2005. It still remains a question if one could
prove (\ref{35}) using Stein's method, but by another approach, perhaps
by a
non-uniform bound on $\Delta_{ik} g_A$ or by a different Stein equation.

\subsection{Other generalisations}\label{sec7.4}

There are two interesting generalisations of Poisson approximation
which we
will not discuss in this paper but will mention in passing. First,
\citet{Brown2001} developed probabilistic methods for approximating general
distributions on non-negative integers with a new family of distributions
called polynomial birth-death distributions. These distributions
include as
special cases the Poisson, negative binomial, binomial and hyper-geometric
distributions. Second, \citet{Peccati2011} combined Stein's method with the
Malliavin calculus of variations to study Poisson approximation for
functionals of general Poisson random measures. This is a follow-up to his
very successful work (see \citet{Nourdin2012}) in normal approximation
for Gaussian functionals using Stein's method and the Malliavin
calculus. Both
the work of \citet{Brown2001} and of \citet{Peccati2011} open up new
domains for Poisson-related approximations and applications of Stein's method.

\section*{Acknowledgments}

We would like to thank Aihua Xia for helpful discussions on generalisations
and open problems and Andrew Barbour for his very helpful comments.



\printhistory

\end{document}